\documentclass[11pt]{article}
\newcommand{\proof}{\noindent {\bf Proof: }}

\newtheorem{statement}{Statement}
\newtheorem{defi}{Definition}
\def\qed{\hfill $\Box$}

\usepackage{amssymb}

\begin{document}
\title{Projection pencil of quadrics and Ivory theorem}
\author{\'A.G.Horv\'ath\\ Department of Geometry, \\
Budapest University of Technology and Economics,\\
H-1521 Budapest,\\
Hungary}
\date{May 15, 2009}

\maketitle

\begin{abstract}

We rewrite the property of confocality with respect to a pseudo-Euclidean space. Our observation is that the generalized definition of confocality in $\cite{stachel}$ does not give back to the original definition of confocality of Euclidean conics. The "confocality property" of $\cite{stachel}$ can be get from our definition in the case when the projection transformation is regular, it is the identity transformation of the space. Our examination concentrate to the singular cases, too because they are playing important role in the investigations of essential examples.
\end{abstract}



\section{Introduction}

Two hundred years ago J.Ivory published a geometric theorem (\cite{ivory}) which importance can be measured on its influence to physic. Ivory theorem states that the level surfaces of the gravitational potential in the exterior domain of an elliptic layer (which is an infinitely thin layer between two similar concentric ellipsoids) are confocal ellipsoids. In mechanics the pair of this theorem the so-called Newton theorem (generalizing the theorem on attraction of spheres to elliptic layers) which says that the gravitational potential inside the elliptic layer is constant.

Some important moment of the long history of this theorem can be red in \cite{stachel}. We would like to mention here two omitted papers: the first one is a characterization of confocal conics in pseudo-Euclidean space due to G.Birkhoff and R.Morris in \cite{birkhoff}, and the second one written by V.V.Kozlov giving analogs of the mechanical variation of Newton and Ivory theorems in spaces of constant curvature of dimension 3. (\cite{kozlov})

In this paper we will rewrite the property of confocality. The origin of this investigation is the concept of confocality introduced in \cite{stachel}. With respect to this definition it was proved that the classical theorem of Ivory is true.
Our observation is that the introduced definition of confocality does not give back to the original definition of confocality of Euclidean conics, but most statements of \cite{stachel} are true for all class of pencil of quadrics called by me "projection pencil of quadrics". We also note that without projective apparats H.Stachel immediately proved Ivory's theorem in the physical Minkowski plane in \cite{stachel1} for standard confocal conics as we would like to also generalize in this paper. So our results can be considered also the common extraction of the results of the papers \cite{stachel1} and \cite{stachel}. In my opinion there is no uniquely determined natural generalization of "confocality property" thus we suggest here such definition which contains all the well-known situations as special cases.

Our paper basically follows the building up of the paper $\cite{stachel}$ uses its terminology, statements and proofs. We prove only such results which are not immediate consequence of the previously proved ones.

\subsection{Notation and Terminology}

We use some concepts without definition as the concepts of real vector spaces, direct sum of subspaces, linear and bilinear mapping,
quadratic forms, projective space, kernel and rank of a linear mapping, respectively.

\begin{description}

\item[$V$, $V^\star$, $L(V)$:] An vector space, its dual space and the space of linear transformations of $V$.
\item[$\dim(V)$]: The dimension of the vector space $V$ in this paper is equal to $(n+1)$.
\item[$P(V)$, $P(L(V))$:] The projective space defined on $V$ and $L(V)$, respectively. In our paper $\dim (P(V))=n$.
\item[$q$,$Q(V)$, $P(Q(V))$:] A quadratic form and the vector resp. projective vector spaces of all quadratic forms of $V$.
\item[$\mbox{ Ker }L$, $\mbox{ Im }L$:] The kernel and image space of the linear transformation $L$.
\item[$\mathbb{C}$, $\mathbb{R}$, $\mathbb{R}^{n}$:] The complex line, the real line and the $n$-dimensional
    Euclidean space, respectively.
\item[{$\langle\cdot,\cdot\rangle$}:] The notion of indefinite inner (scalar) product.
\item[$\Phi(x,y)$,$\Phi$, $\Omega$:] The notation of a symmetric bilinear function, its zero set called by the corresponding quadric and the absolute quadric associated to a fixed projection with $n$ or $(n+1)$-dimensional image space, respectively.
\item[$l,g,id$, $L,G,E$:] The notation of linear mapping $l,g,id:V\longrightarrow V$ (and also their projective classes in $P(L(V))$) and their $(n+1)\times (n+1)$ matrices with respect to a homogeneous coordinate system, respectively.
\item[$P$:] Is a projection transformation of $V$. It holds (and defined by holding) the equality $P^2=P$.
\end{description}

\section{Quadrics in a finite dimensional projected space}

There is a well-known building up of projective quadrics in an $n$-dimensional projective space $P(V)$  over a commutative field $K$ of characteristic different from 2. (See \cite{berger}) If $V$ is a vector space over the field $K$ we can consider the elements of the set $Q(V)$ of homogeneous polynomial of degree $2$ as quadratic forms. By definition for such a $q:V\longrightarrow K$ there is a symmetric, bilinear map $\Phi$ on $V\times V$ into $K$ which is the polar form of $f$, holding the equality:
$$
 \Phi(x,y)=\frac{1}{2}(f(x+y)-f(x)-f(y))
$$
for all $x,y\in V$. Denote by $P(V)$ the $n$-dimensional projective space associated to the $n+1$-dimensional vector space $V$ and the canonical projection $p:V\setminus 0\longrightarrow P(V)$, respectively. The isotropic cone of $q$ defined by $p(q^{-1}(0)\setminus 0)$ does not change when we replace $q$ by $kq$ $k\in K^\ast$. This means that such a projection is actually associated to a point of the projective space $P(Q(V))$. The projective quadrics (in algebraic sense) are the elements of $P(Q(V))$. A quadric is a conic when $n=2$. It is proper if it has a non-degenerate equation; otherwise it is called degenerate. The classification of quadrics in real and complex cases are known, we have to determine the orbits of $P(Q(V))$ under the action of the isometry group $GP(V)$. (This is the projective group of $P(V)$ acts on $P(Q(V))$. If we consider that the fields of coefficients is $\mathbb{C}$ then we have exactly $(n+1)$ orbits, classified by the rank $k$, where $1\leq k\leq (n+1)$, and in the real case the orbits are classified by pairs $(r,s)$ such that $1\leq s\leq r\leq n+1$. In particular, there exists $\frac{1}{2}n(n+1)+1$ types of proper quadrics. (See 14.1.5.1. Theorem in \cite{berger}.)

In geometric point of view we can see the quadrics as the zero set of quadratic forms (or the zero set of its symmetric bilinear forms). Fixing a regular symmetric bilinear form as an indefinite inner product $\langle \cdot,\cdot \rangle$ all quadrics can be regarded as the zero set of a symmetric bilinear function $\langle x, l(y) \rangle$ where $l$ is a self-adjoint transformation with respect to the product $\langle \cdot,\cdot \rangle$. Every such transformation represent a quadric but the elements of the projective class of $l$ belong to the same quadric. It is also possible that the same quadric associated to two non-equivalent self-adjoint transformation. Every linear transformation $l$ of $V$ means an endomorphism of the corresponding projective space $P(V)$, too. We do not use another notation for the endomorphisms and the linear transformations corresponding to them. In a numerous calculation we can use matrix representation of the transformations changing suitable basis for it. There is a basic result on the characterization of the fixed symmetric bilinear form $\langle \cdot,\cdot \rangle$ and the self-adjoint linear mapping $l$ by simultaneous normal form - meaning the choice of a basis such that both have simple coordinate elements with respect to this basis. (See \cite{gohberg1} Th. 5.3. or \cite{gohberg}.) The matrix representation of linear transformations denoted by capital letters.

Now we can introduce the dual of a quadric which is a quadric of the dual space.
Two points $x$ and $y$ of $P(V)$ are conjugate with respect to the symmetric bilinear form $\Phi$ if $\Phi(x,y)=0$. The set of points conjugate to $x$ is a subspace of $P(V)$ of dimension at least $(n-1)$.  Since the $(n-1)$-dimensional subspaces in $P(V)$ are the zero sets of linear forms $a^\ast \in P(V^{\ast})$, the set of conjugate points is a point $x^\ast$ of $P(V^{\ast})$. The bilinear form
$$
\widehat{\Phi}(x^\ast,y^\ast):=\Phi(x,y)
$$
is called the dual of the form $\Phi$ and the corresponding quadric of $P(V^\ast)$ is the dual of the quadric defined by $\Phi$. If there is a  correspondence (distinct for the duality) between the dual form and a bilinear symmetric quadratic form of $P(V)$ then we say that we gave a pullback of the dual form into $P(V)$ We denote by $\widetilde{\Phi}$ the pullback of $\widehat{\Phi}$. To see this concretely assume that the bilinear form is $\Phi (x,y):=x^TTy$ by a symmetric matrix $T$ and its dual quadric is the zero set of $\widehat{\Phi }(x^\ast,y^\ast)$. To determine it we have to give the duality map.
Let $(\cdot)^\ast:V\longrightarrow V^\ast$ defined by
$$
x\mapsto x^\ast =\Phi_{x}(\cdot):=\Phi (x,\cdot).
$$
If $L$ is a linear transformation of $V$ let denote by $L^T$ the transposed of $L$ with respect to the fixed product. Then
$$
\Phi(L(x),L(y))=(L(x))^TTL(y)=x^T(L^TTL)y
$$
and
$$
L(x)^\ast=\Phi_{L(x)}(\cdot)=\Phi (L(x),\cdot).
$$
If $u\in \mbox{ Im }T$ then there is an $x\in V$ for which $T(x)=u$ and we can use the notation $T^{-1}u:=x+\mbox{ Ker }T$. Thus for $u,v\in \mbox{ Im }T$ we can define a pullback $\widetilde{\Phi}$ of the dual quadric $\widehat{\Phi}$. In fact,
$$
\widehat{\Phi} (x^\ast,y^\ast)=\widehat{\Phi} ((T^{-1}(u))^\ast,(T^{-1}(v))^\ast)=\widehat{\Phi} (\Phi_{T^{-1}u}(\cdot ),\Phi_{T^{-1}v}(\cdot ))=u^T((T^{-1})^TT(T^{-1}))v,
$$
thus the definition
$$
\widetilde{\Phi}(u,v):=u^T((T^{-1})^TT(T^{-1}))v
$$
gives a pullback of the dual form. The set valued mapping $T^{-1}$ gives an isomorphism on $\mbox{ Im }T$ to the factor space $V/\mbox{ Ker }T$, thus the mapping $T(T^{-1})$ is identity on $\mbox{ Im }T$. So the pullback of the dual quadratic form using the symmetricity of $T^{-1}$ can be defined on Im $T$ by
$$
\widetilde{\Phi}(u,v):=u^T(T^{-1})^Tv=u^TT^{-1}v.
$$
We extract the linear transformation $T^{-1}:\mbox{ Im }\longrightarrow V/\mbox{ Ker }T$ to $V$ by the zero map of the complementary subspace of $\mbox{ Im }T$ giving  those singular transformation which define the singular symmetric bilinear form $\widetilde{\Phi}(u,v)$ of $V$.

\section{The confocality of conics}

In this section let the dimension $n$ is equal to two. To define confocality in the general case, we shortly recall the description of Euclidean confocal conics. Working in Euclidean homogeneous coordinates, finite points and asymptotic directions ("points at infinity") are and given by column vectors. A row vector specify as a line with normal $n$. The line at infinity contains the infinite points. Change-of-basis transformations are matrices acting by left multiplication on points and by right multiplication by the inverse on planes so that point-line products are preserved. Euclidean transformations take the form `
$$
T=\left( \begin{array}{cc}
          A& b\\
         0& 1
        \end{array}
 \right)
$$ where $A$ is a rotation matrix $b$ a translation vector. $A$ becomes a re-scaled rotation for scaled Euclidean or similarity
transformations, and an arbitrary nonsingular $3\times 3$ matrix for affine ones.
Consider first the absolute quadric corresponds to the symmetric rank $2$ matrix
$$
\Omega=\left( \begin{array}{ccc}
          1&0& 0\\
          0&1&  0\\
         0&0& 0
        \end{array}
 \right).
$$
The coordinates of a point $x$ of the absolute holds the equality $x_1^2+x_2^2=0$ thus it has only one point $(0,0,1)^T$. The pullback of the dual form is
$$
\widetilde{\Phi}(u,v)=(u_1,u_2)\left( \begin{array}{ccc}
          1&0&0\\
          0&1& 0\\
          0&0& 0
        \end{array}
 \right)\left( \begin{array}{c}
          v_1 \\
         v_2
        \end{array}\right),
$$
where $u=(u_1,u_2,0)^T$ and $v=(v_1,v_2,0)^T$ are elements of $\mbox{ Im }\Omega$. The pullback of the dual quadric is also the absolute, so it has only one element associated to the projective point $(0,0,1)^T$. Of course, the $3\times 3$ identity matrix as a selfadjoint transformation identical with its inverse showing that it defines a self-dual conics which is the empty set.
It is easy to see that confocal conics with foci $(0,\pm c)$ can be written in a suitable homogeneous coordinate system of the form:
$$
0=(x_1,x_2,x_3)\left( \begin{array}{ccc}
          \frac{1}{c^2+\lambda }& 0 & 0\\
         0& \frac{1}{\lambda } & 0\\
         0& 0 & -1
        \end{array}\right)\left(\begin{array} {c}
                          x_1\\
                          x_2\\
                          x_3
\end{array}\right).
$$
Since this quadratic form defined by a regular matrix the pullback of its dual is:
$$
0=(x_1,x_2,x_3)\left( \begin{array}{ccc}
          c^2+\lambda & 0 & 0\\
         0& \lambda  & 0\\
         0& 0 & -1
        \end{array}\right)\left(\begin{array} {c}
                          x_1\\
                          x_2\\
                          x_3
\end{array}\right),
$$
showing that the pullback of the dual forms of confocal conics giving a linear subset of all forms, containing the pullback of the dual of the absolute. In fact,
$$
\alpha\left( \begin{array}{ccc}
          c^2+\lambda_1 & 0 & 0\\
         0& \lambda_1  & 0\\
         0& 0 & -1
        \end{array}\right)+
 \beta  \left( \begin{array}{ccc}
          c^2+\lambda_2 & 0 & 0\\
         0& \lambda_2  & 0\\
         0& 0 & -1
        \end{array}\right)=
$$
$$
       = \left( \begin{array}{ccc}
          (\alpha +\beta)c^2+(\alpha\lambda_1+\beta\lambda_2) & 0 & 0\\
         0& (\alpha\lambda_1+\beta\lambda_2) & 0\\
         0& 0 & -(\alpha +\beta)
         \end{array}\right)
         =
$$
$$
 =\left \{
         \begin{array}{lcl}\left( \begin{array}{ccc}
          c^2+\lambda_3 & 0 & 0\\
         0& \lambda_3  & 0\\
         0& 0 & -1
                \end{array}\right)& \mbox{ if }&  (\alpha +\beta)\neq 0\\
         \left( \begin{array}{ccc}
          \lambda_3 & 0 & 0\\
         0& \lambda_3  & 0\\
         0& 0 & 0
         \end{array}\right) & \mbox{ if } & (\alpha +\beta)=0.
         \end{array}\right.
$$
On the other hand this linear hull does not contain the $3\times 3$ identity transformation. This implies that the definition of \cite{stachel} is not a generalization of the usually confocality of conics. As it can be seen easily, for a non-zero $c$ the connection between the distances $2c_\alpha $ and $2c$ of foci of the conics $\Phi_\alpha$ and $\Phi$ for which $\Phi_\alpha ^{-1}=\alpha \Phi +\beta I$ by the $3\times 3$ identity matrix $I$ is
$$
c_\alpha ^2=\frac{\alpha}{\alpha-\beta}c^2
$$
showing that with respect to the pencil of conics of paper \cite{stachel} the distances of foci are not constant.

Our second example is the confocal family of conics with two axes of symmetry belonging to the pseudo-Euclidean (Minkowski) plane. In the paper \cite{birkhoff} it was called by relativistic confocal conics in space time, and shown the relativistic conics are geometrically tangent to the null lines (isotropic lines) through the foci. In paper \cite{stachel1} it was type B, and we can find two nice figure on it. Using our setting up we can see the followings:

The regular bilinear function in a cartesian homogeneous coordinate system of the embedding Euclidean plane is:
$$
\langle x,y \rangle=(x_1,x_2,x_3)\left(\begin{array}{ccc}
                                                                      1& 0& 0\\
                                                                      0&-1& 0\\
                                                                      0&0&1
                                                                      \end{array}\right)\left(\begin{array}{c}
                                                                      y_1\\
                                                                      y_2\\
                                                                      y_3
                                                                      \end{array}\right).
$$
The absolute $\Omega$ can be considered as the zero set of the bilinear function:
$$
\langle (x_1,x_2,x_3),\left(\begin{array}{ccc}
                        1& 0& 0\\
                        0&1& 0\\
                        0&0&0
                                 \end{array}\right)\left(\begin{array}{c}
                                                                      y_1\\
                                                                      y_2\\
                                                                      y_3
                                                                      \end{array}\right)\rangle=x_1y_1-x_2y_2,
$$
defined by the projection transformation
$$
P=\left(\begin{array}{ccc}
                        1& 0& 0\\
                        0&1& 0\\
                        0&0&0
                                 \end{array}\right).
$$
Since we can write the equation:
$$
\frac{x_1^2}{\sigma}+\frac{x_2^2}{\tau }=1 \mbox{ with } \sigma \tau (\sigma +\tau),
$$
by the selfadjoint transformation
$$
G=\left(\begin{array}{ccc}
                        \frac{1}{\sigma}& 0& 0\\
                        0&-\frac{1}{\tau}& 0\\
                        0&0&-1
                                 \end{array}\right),
$$
into the form
$$
\langle x,Gx\rangle=0.
$$
The normal form of pencil of the corresponding confocal conics can be got from the equality
$$
0=\langle x,(G^{-1}-tP)^{-1}x\rangle=(x_1,x_2,x_3)\left(\begin{array}{ccc}
                        \frac{1}{\sigma-t}& 0& 0\\
                        0&\frac{1}{-(-\tau-t)}& 0\\
                        0&0&-1
                                 \end{array}\right)\left(\begin{array}{c}
                                                                      x_1\\
                                                                      x_2\\
                                                                      x_3
                                                                      \end{array}\right).
$$
With respect to the original inhomogeneous cartesian coordinates it is:
$$
\frac{x_1^2}{\sigma-t}+\frac{x_2^2}{\tau +t}=1 \mbox{ for } t\in \mathbb{R} \setminus \{\sigma ,\tau\}.
$$

The elliptic and hyperbolic case can be considered parallely. We consider special ovals of the projective space as an intersection of a family of quadratic cones and the model planes, respectively. In the first case, the point of the model of the elliptic plane are the opposite pairs of the unit sphere of $V$ where the lengthes of the vectors calculated by the bilinear function
$$
\langle x,y \rangle=(x_1,x_2,x_3)\left(\begin{array}{ccc}
                                                                      1& 0& 0\\
                                                                      0&1& 0\\
                                                                      0&0&1
                                                                      \end{array}\right)\left(\begin{array}{c}
                                                                      y_1\\
                                                                      y_2\\
                                                                      y_3
                                                                      \end{array}\right).
$$
In the second case the bilinear function is
$$
\langle x,y \rangle=(x_1,x_2,x_3)\left(\begin{array}{ccc}
                                                                      1& 0& 0\\
                                                                      0&1& 0\\
                                                                      0&0&-1
                                                                      \end{array}\right)\left(\begin{array}{c}
                                                                      y_1\\
                                                                      y_2\\
                                                                      y_3
                                                                      \end{array}\right),
$$
and the points of the model is the opposite pairs of points of the hyperboloid containing the vectors with imaginary unit lengthes. The family of cones defined by the equalities:
$$
\frac{x_1^2}{c^2}+ \frac{x_2^2}{c^2-\beta^2}\pm \frac{x_3^2}{c^2+\gamma^2}=0,
$$
where $c$ is a parameter and $\beta^2\pm \gamma^2=\pm 1$ in the two respective cases. The definition in the elliptic case is also motivated by physical argument of a gravitating arc, as we can see in \cite{kozlov}.

By the selfadjoint transformation
$$
G=\left(\begin{array}{ccc}
                        \frac{1}{c^2}& 0& 0\\
                        0&\frac{1}{c^2-\beta^2}& 0\\
                        0&0&\frac{1}{c^2+\gamma^2}
                                 \end{array}\right),
$$
we can write
$$
\langle x,Gx\rangle=0.
$$
The normal form of pencil of the corresponding confocal conics can be got from the equality
$$
0=\langle x,(G^{-1}-tI)^{-1}x\rangle=(x_1,x_2,x_3)\left(\begin{array}{ccc}
                        \frac{1}{c^2-t}& 0& 0\\
                        0&\frac{1}{(c^2-t)-\beta^2}& 0\\
                        0&0&\frac{1}{(c^2-t)+\gamma^2}
                                 \end{array}\right)\left(\begin{array}{c}
                                                                      x_1\\
                                                                      x_2\\
                                                                      x_3
                                                                      \end{array}\right),
$$
by the matrix
$$
I=\left(\begin{array}{ccc}
                                                                      1& 0& 0\\
                                                                      0&1& 0\\
                                                                      0&0&1
                                                                      \end{array}\right)
$$
of the identical transformation. Thus the normal forms of these pencil of conics are
$$
\frac{x_1^2}{c^2-t}+ \frac{x_2^2}{(c^2-t)-\beta^2}\pm \frac{x_3^2}{(c^2-t)+\gamma^2}=0 \mbox{ where } x_1^2+x_2^2\pm x_3^2=\pm 1.
$$

We can conclude that for the description of confocality we have to use both of the singular and nonsingular projection transformation of the space $V$. This motivates our further examination.

On the Euclidean and hyperbolic plane the conics of different types are the ellipses (there is no ideal points) and hyperboles (with two ideal points), respectively. The singular quadrics (determining the common line of the two foci) divided the all family into these two types. This situation can be observed in the elliptic case, too, but there is no other possibility to distinguish the getting ovals to each other.

In the previously investigated case of the pseudo-Euclidean plane there are three types of conics as it can be seen either in \cite{stachel1} or in \cite{birkhoff}.

\section{Ivory property, projection pencil of quadrics and $p$-quadrics}

The planar Euclidean version of Ivory's theorem states that the two diagonals of any curvilinear quadrangle formed by four confocal conics have the same length. With respect to the scalar product the equality means:
$$
\rho^2 (x,y')=\langle x-y',x-y'\rangle =\langle y-x',y-x'\rangle =\rho^2 (x',y),
$$
where the pairs of points $\{x,y\},\{x,x'\},\{y,y'\},\{x',y'\}$ are on four confocal quadrics which intersect to each other in the examined points.
Also in Euclidean space there is an equivalent reformulation of this theorem on the language of affine mapping, since if we have two confocal conics of the same type (e.g. ellipses) then there exists an affine mapping $l$ with the property, that whenever a conic of other type (hyperbole) intersects the first ellipses in a point then it intersects the other one in those point which is the image of the first one by this mapping. Both intersection are orthogonal and now Ivory theorem states that:
$$
\rho(x,l(y))=\rho(y,l(x)),
$$
where $\rho $ is the Euclidean distance. In a projective pseudo-Euclidean space the hyperbolic and elliptic metric (based on the inner product) is a function of the quantity
$$
\rho(x,y):=\frac{\langle x,y\rangle}{\sqrt{|\langle x,x\rangle\langle y,y\rangle}|},
$$
so it was a natural conception  to substitute the lengthes of the examined diagonals by $\rho $. (See in \cite{stachel}.)

We call by projection transformation (briefly projection) a linear transformation $p:V\longrightarrow V$ if it holds the equality $p^2=p$. The vector space $V$ can be regarded as the direct sum of its subspaces $\mbox{ Ker }p$ and $\mbox{ Im }p$. The projection transformation restricted to its image space is the identity one, and to its kernel is the zero mapping, respectively.

First we note that the property of a linear transformation named by projection is not a "projective property", in the projective class of the projection $p$ the only projection transformation is $p$. In fact, if $\lambda\neq 0,1$ then $(\lambda p)^2=\lambda^2 p\neq \lambda p$, showing that the transformation $\lambda p$ is not a projection. We say that a point of $P(L(V))$ is a projection if there is a representant of its class which is a projection.

\begin{defi}
A maximal set of quadrics called by pencil of quadrics if the selfadjoint linear transformations corresponding their duals belong a two-dimensional subspace of the vector space $L(V)$. We say that the pencil of quadrics is a projection pencil associated to the projection $p$, if the corresponding two-space contains $p$. Within a family of projection pencil of quadrics spanned by the invertible linear transformation $l_0$ and the projection $p$, the connected components of $\{(\lambda ,\nu)\mbox{ } |\mbox{ } \lambda (l_0^{-1}+\mu p)\}$ correspond to quadrics of different type.
\end{defi}

Lemma 9 and 10 in \cite{stachel} give a representation theorem showing that the special case investigated in lemmas 6,7 and 8 is actually the general case to need for Ivory's theorem. Our purpose to describe the general situation when this representation is possible. For this reason we introduce a new concept, the concept of $p$-quadric.

\begin{defi}
Let $p$ be a projection. The quadric $\Phi $ generated by the selfadjoint transformation $g$ is a $p$-quadric  if for every $w=u+v\in V$ for which $u\in \mbox{ \textrm{Im} }p$ and $v\in \mbox{ \textrm{Ker} }p$ we have $g(u+v)=pg(u)-v$.
\end{defi}

We remark that every quadric is an $id$-quadric and if $\Phi$ is a $p$-quadric then $\mbox{ Im }p$ is an invariant subspace of $g$. In fact, if $u\in \mbox{ Im }p$ then
$$
g(u)=p(g(u)),
$$
proving that $g(u)$ is also an element of $\mbox{ Im }p$.

By the method of H.Stachel and J.Wallner we can prove the following theorem.

\noindent {\bf Theorem}(A generalization of Ivory's theorem){\bf :} {\em Let $P(V)$ be a projective space  with metric
$$
\delta (x,y):=\frac{\langle x,y\rangle}{\sqrt{|\langle x,x\rangle\langle y,y\rangle|}}
$$
where $\langle\cdot,\cdot\rangle$ is a fixed indefinite inner product of $V$ and let $p$ be a projection of $V$. Denoted by $\Phi_0=\{ x \mbox{ }|\mbox{ } 0=\langle x, l(x)\rangle\}$ and $\Phi_1=l_1(\Phi_0)$ two regular $p$-quadric of the same type, belonging to the projection pencil associated to $l$ and $p$. Then there is a smooth family $\Phi_\lambda =l_\lambda (\Phi_0)$ ($0\leq \lambda \leq 1$ of $p$-quadrics of this pencil, such that $l_\lambda $ is selfadjoint and has the Ivory property:
$$
\delta(x,l_\lambda (y))=\delta (l_\lambda (x),y) \mbox{ for all } x,y\in \Phi_0\cap \mbox{ Im }p.
$$
Any further $p$-quadric $\Psi $ corresponding to the same projection pencil and containing a point $x\in \Phi_0\cap \mbox{ \textrm{Im} }p$, also contains the entire path $l_\lambda (x)$, which intersects all quadrics $\Phi_\lambda $ orthogonally in $\mbox{ \textrm{Im} }p$.}

\section{The complete list of the cited definitions and statements}

Our paper to be self-contained we give the complete list and original numeration of statements and definitions in \cite{stachel}:

\noindent{\bf Definition 1.} {\em A (nondegenerate) quadric $\Phi$ is the zero set of a (nondegenerate) symmetric bilinear form $\sigma (x,y)=\langle x, l(y) \rangle$, with a selfadjoint (nonsingular) linear endomorphism $l$. $l=id$ corresponds to the absolute quadric $\Omega $, the set of absolute points.}

\noindent{\bf Definition 2.} {\em The quadric $\widehat{\Phi}=\{v \mbox{ }|\mbox{ } \widehat{\sigma}(v,v)=\langle v,l^{-1}(v)\rangle=0 \}$ in the dual space is called the dual of the original quadric $\Phi $ defined by $\sigma (v,w)=\langle v,l(w)\rangle$.}

\noindent{\bf Definition 3.} {\em If $k$ is a linear endomorphism and the quadric $\Phi $ is given by the endomorphism $l$, then we define the dual $k$-image of $\Phi $ to have the equation
$$
\widehat{\sigma}(v,v)=0, \mbox{ with } \widehat{\sigma}(v,v)=\langle v,kl^{-1}k^{\ast}(w)\rangle
$$
$\widehat{\sigma}$ is understood to apply to gradients.}

\noindent{\bf Definition 4.} {\em $\Phi _0$ and $\Phi _1$ are said to be confocal (or homofocal ), if one of the following equivalent conditions holds true:

\noindent (i) the bilinear forms $\widehat{\sigma_0}$, $\widehat{\sigma}_1$, $\widehat{\langle \cdot,\cdot \rangle}=\langle \cdot,\cdot \rangle$ are linearly dependent,

\noindent (ii) the linear endomorphisms $l^{-1}_0$, $l^{-1}_1$, $id$ are linearly dependent,

\noindent (iii) the coordinate matrices $Q^{-1}_1$, $Q^{-1}_2$, $H^{-1}$ are linearly dependent.

\noindent The family of quadrics $\Phi $ confocal to $\Phi_0$ is defined by endomorphisms $l$ which satisfy $l^{-1}=\lambda l^{-1}_0+\mu id$, $(\lambda, \mu)\in \mathbb{R}^{2}$, $\lambda \neq 0$.}

\noindent{\bf Definition 5.} {\em Within the family of confocal bilinear forms spanned by $l_0$, the connected components of $\{(\lambda, \mu) \mbox{ }|\mbox{ } \lambda(l^{-1}_0+\mu id) \mbox{ nonsingular}\}$ correspond to quadrics of different types.}

\noindent{\bf Lemma 3} {\em In the $n$-dimensional elliptic or hyperbolic space $(n>1)$ all confocal families possess at least two types of quadrics.}

\noindent{\bf Lemma 4} {\em If confocal quadrics $\Phi_0$ and $\Phi_{\lambda}$ intersects, they do so orthogonally.}

\noindent{\bf Lemma 5} {\em Assume that $\Phi $ is a quadric, possible singular but not contained in a hyperplane, and that there is a mapping $x\longrightarrow x'$ such that
$$
\langle x'_1,x_2\rangle=\langle x_1,x'_2\rangle \mbox{ for all } x_1,x_2\in \Phi,
$$
then there is a selfadjoint linear endomorphism $l$ of $\mathbb{R}^{n+1}$ such that $x'=l(x)$ for all $x\in \Phi$.}

\noindent{\bf Lemma 6} {\em If the linear endomorphism $l$ is selfadjoint, then the quadric
$$
\Phi_0:\sigma (x,x):=\langle x,x \rangle-\langle l(x), l(x) \rangle =0
$$
together with its $l$-image $\Phi_1$ has the Ivory property
$$
\delta(l(x),y)=\delta(x,l(y)) \mbox{ for all } x,y\in \Phi \mbox{ with } \langle x,x \rangle,\langle y,y \rangle\neq 0.
$$
The restriction of $l$ to any linear subspace contained in $\Phi_0$ is isometric in the sense of $\delta $.}

\noindent{\bf Lemma 7} {\em Assume that $l$ is selfadjoint and that the quadric $\Phi_0$ given in Lemma 6 is regular. Then $\Phi_0$ and $\Phi_1=l(\Phi_0)$ are confocal. (The dual of $l(\Phi_0)$ defined by the endomorphism $lg_0^{-1}l^{\ast}$ if $\Phi_0$ given by $g_0$.) }

\noindent{\bf Lemma 8} {\em If $l$ is selfadjoint, then in most cases the quadric $\Phi_0$ as defined in Lemma 6 is of the same type as $l(\Phi_0)$ provided both are regular. Different types are only possible when the normal form of $l$ contains a block matrix $R_2(0,b)$ or $R_{2k}(0,b,1)$.}

\noindent{\bf Lemma 9} {\em Consider two regular confocal quadrics $\Phi_0$, $\Phi_1$ which are of the same type. Then there is a selfadjoint endomorphism $l$ such that $\Phi_1=l(\Phi_0)$ and the equation of $\Phi_0$ is given by $\langle x,x \rangle-\langle l(x), l(x) \rangle =0$.}

\noindent{\bf Lemma 10} {\em We use the notation of the proof of Lemma 9. There is $\delta > 0$ such that $id- \lambda g_0$ has a square root which smoothly depends on $\lambda $, for $-\delta<\lambda <1+\delta $.}

\noindent{\bf Lemma 11} {\em Suppose that $P,\Phi_0,\Phi_1,g_0,g_1,l$ are as in Lemma 9 and its proof. Then there is a smooth family $l_\lambda$ of transformations with $l_0=id$ and $l_1=l$, such that the quadric $\Phi_\lambda=l_\lambda(\Phi_0)$ is defined by the endomorphism $g_\lambda $ with
$$
g_\lambda ^{-1}=g_0^{-1}-\lambda id.
$$
All quadrics $\Phi_\lambda$ are confocal with $\Phi_0$. They orthogonally intersect the path $l_\lambda(x)$ of a point $x\in \Phi_0$.}

\noindent{\bf Lemma 12} {\em We use the notations of Lemma 11 and consider the quadrics $\Phi_\lambda $, defined by endomorphisms $g_{\lambda}$. If $\Psi \neq \Phi_0$ is confocal with $\Phi_0$, and $x\in \Phi_0\cap \Psi$, then also $l_{\lambda} (x)\in \Psi$.}

\noindent{\bf Theorem 2} {\em Let $P(V)$ be a projective space  with metric
$$
\delta (x,y):=\frac{\langle x,y\rangle}{\sqrt{|\langle x,x\rangle\langle y,y\rangle|}}
$$
where $\langle\cdot,\cdot\rangle$ is a fixed indefinite inner product of $V$. Denoted by $\Phi_0$ and $\Phi_1=l_1(\Phi_0)$ two regular confocal quadrics of the same type. Then there is a smooth family $\Phi_\lambda =l_\lambda (\Phi_0)$ ($0\leq \lambda \leq 1$) of quadrics confocal with $\Phi_0$ and $\Phi_1$, such that $l_\lambda $ is selfadjoint and has the Ivory property:
$$
\delta(x,l_\lambda (y))=\delta (l_\lambda (x),y) \mbox{ for all } x,y\in \Phi_0\cap \mbox{ Im }p.
$$
Any further quadric $\Psi $ confocal with $\Phi_0$ which contains a point $x\in \Phi_0$ contains the entire path $l_\lambda (x)$, which intersects all quadrics $\Phi_\lambda $ orthogonally.}

\section{The proof of the Theorem}

We now modify the statements of the previous section if it need. Since a projection transformation $p$ restricting to $\mbox{ Im }p$ is the identity $id$ , the proof of Lemma 3 is valid. Lemma 4 for our projection pencil of quadrics can be formulated in the following way:

\begin{statement}[Lemma 4']
If quadrics $\Phi_0$ and $\Phi_{\lambda}$ corresponding to a projection pencil intersect, they do so orthogonal with respect to the quadratic form of $p$. More precisely, if $x$ is a common point then we have $0=\langle p(g_0(x)),g_{\lambda}(x)\rangle  $.
\end{statement}

\proof In fact, we have
$$
0=\langle x,g_0(x)\rangle  =\langle x,v\rangle  =\langle g_{\lambda}^{-1}(w),v\rangle  =\langle (g_{0}^{-1}+\mu p)(w),v\rangle  =
$$
$$
=\langle g_{0}^{-1}(w),v\rangle  +\mu \langle p(w),v\rangle  =\langle w,g_{0}^{-1}v\rangle  +\mu \langle p(w),v\rangle  =
$$
$$
=\langle g_{\lambda}(x),x\rangle  +\mu \langle p(w),v\rangle  =\mu \langle p(w),v\rangle  ,
$$
as we stated.

\qed

Lemma 5 on Ivory property is valid, again. In Lemma 6, Lemma 7 and Lemma 8 we change the selfadjoint transformation $l$ to the selfadjoint transformation $l'=lp+(id-p)$ and we consider the quadrics $\Phi_0'$ with equation:
$$
\langle p(x),p(x)\rangle  -\langle l'(x),l'(x)\rangle  =0.
$$
We remark that on $\mbox{ Im }p\cap \Phi_0'=\mbox{ Im }p\cap \Phi_0$ where $\Phi_0$ defined by the equality
$$
\langle x,x\rangle  -\langle l(x),l(x)\rangle  =0
$$
and  $\mbox{ Ker }p\cap \Phi_0'=\{x\in \mbox{ Ker }p \mbox{ }|\mbox{ } \langle x,x\rangle  =0\}$.

As it can be seen easily the following variation of Lemma 6 is true for every projection pencil of quadrics:

\begin{statement}[Lemma 6']
If the linear endomorphism $l$ is selfadjoint and invariant on the subspace $\mbox{ Im }p$ then the quadric
$$
\langle p(x),p(x)\rangle  -\langle l'(x),l'(x)\rangle  =0
$$
together with its $l'$ image $\Phi_1=l'(\Phi_0)$ has the Ivory property
$$
\delta(l'(x),y)=\delta(x,l'(y)) \mbox{ for all } x,y\in \Phi \mbox{ with } \langle x,x\rangle  , \langle y,y\rangle   \neq 0.
$$
The restriction of $l'$ to any linear subspace contained in $\Phi_0$ is isometric in the sense of $\delta $.
\end{statement}

The following modification is more interesting:

\begin{statement}[Lemma 7']
Assume that $\mbox{ Im }p$ is an invariant subspace of $l$, and that the quadric $\Phi_0$ given by the equality:
$$
\langle p(x),p(x)\rangle  -\langle l'(x),l'(x)\rangle  =0, \mbox{ where } l'=lp+(id-p)
$$
is regular. Then $\Phi_1=l'(\Phi_0)$ is in the projection pencil of $\Phi_0$ and $p$.
\end{statement}

\proof Rewriting the equation of $\Phi_0$ we get:
$$
0=\langle x,p(x)\rangle  -\langle l'(x),l'(x)\rangle  =\langle x,(p-(l')^2)x\rangle  .
$$
$V$ is a direct sum of $\mbox{ Ker }p$ and $\mbox{ Im }p$ for arbitrary $p$. Furthermore $l'=l$ on $\mbox{ Im }p$ and $l'=id$ on $\mbox{ Ker }p$.
The dual of $\Phi_0$ is represented by
$$
0=\langle x,(p-(l')^2)^{-1}x\rangle  ,
$$
and the dual $l'$-image of $\Phi_0$ is according the Def. 3 in \cite{stachel}, defined by
$$
0=\langle x,l'(p-(l')^2)^{-1}l'x\rangle  .
$$
Consider now the transformation
$$
(p-(l')^2)^{-1}-l'(p-(l')^2)^{-1}l'.
$$
Observe that $\mbox{ Im }p$ is an invariant subspace of $(p-(l')^2)$. In fact, for $u\in \mbox{ Im }p$ we also have $(p-(l')^2)(u)=u-l^2(u)\in \mbox{ Im }p$.
Thus $\mbox{ Im }p$ is an invariant subspace of their inverse and for a vector $u\in \mbox{ Im }p$ by the argument of Lemma 7  applying it to the invariant subspace $\mbox{ Im }p$, we have
$$
(p-(l')^2)^{-1}-l'(p-(l')^2)^{-1}l'(u)=((id-l^2)^{-1}-l(id-l^2)^{-1}l)(u)=u.
$$
On the other hand for a vector $v\in \mbox{ Ker }p$
$$(p-(l')^2)(v)=-v
$$
showing $v\in \mbox{ Ker }p$ and $(p-(l')^2)$ is a reflection on $\mbox{ Ker }p$. Thus
$$
((p-(l')^2)^{-1}-l'(p-(l')^2)^{-1}l')(v)=(p-(l')^2)^{-1}(v)-(p-(l')^2)^{-1}(v)=0,
$$
giving the required equality:
$$
(p-(l')^2)^{-1}-l'(p-(l')^2)^{-1}l'=p.
$$
\qed

In further, we consider such selfadjoint transformations only, which leaves invariant $\mbox{ Im }p$.

\begin{statement}[Lemma 8']
If $l$ is selfadjoint with invariant subspace $\mbox{ Im }p$, then in most cases the quadric $\Phi_0$ as defined by $g_0=(p-(l')^2)$ is of the same type as $\Phi_1=l'(\Phi_0)$ provided both are regular. Different types are only possible when the normal form of $l'|_{\mbox{ Im }p}$ contains a block matrix $R_2(0,b)$ or $R_{2k}(0,b,1)$. (See Th.1 in \cite{stachel} or Th.5.3. in \cite{gohberg1}.)
\end{statement}

\proof The convex combination of selfadjoint transformations $g_0^{-1}=(p-(l')^2)^{-1}$ and $g_1^{-1}=l'(p-(l')^2)^{-1}l'$ can be investigated on the way of Lemma 8, using the result of our Statement 2:
$$
g_{\lambda}^{-1}:=(1-\lambda )g_0^{-1}+\lambda g_1^{-1}=g_0^{-1}-\lambda p=(p-(l')^2)^{-1}-\lambda p,
$$
if $0\leq \lambda \leq 1$. For $u\in \mbox{ Im }p$ we have $\left((p-(l')^2)^{-1}-\lambda p\right)(u)=\left((id-l^2)^{-1}-\lambda id\right)(u)$ and the proof of Lemma 8 can be applied. For $v\in \mbox{ Ker }p$ we get that $\left((p-(l')^2)^{-1}-\lambda p\right)(v)=-id(v)$ which is always non-singular.
\qed

Since the quadric $g_0$ using in the previously statements Lemma 6'--Lemma 8' are $p$-quadric we can give representation theorem only $p$-quadrics.

\begin{statement}[Lemma 9']
Consider two regular $p$-quadrics $\Phi_0$ and $\Phi_1$ of a projection pencil which are of the same type with respect to the projection $p$. Then there is a selfadjoint transformation $l$ invariant on the subspace $\mbox{ Im }p$ such that $\Phi_1=l'(\Phi_0)$ with the transformation $l'=lp+(id-p)$ and the equation of $\Phi_0$ is given by
$$
\langle p(x),p(x)\rangle  -\langle l'(x),l'(x)\rangle  =0.
$$
\end{statement}

\proof
Without changing the quadrics we can consider regular representing selfadjoint transformation $g_0$ and $g_1$ for which $\mbox{ Im }p$ is an invariant subspace and
$$
g_0^{-1}-g_1^{-1}=p.
$$
We have to show that there exists $l$ such that it is invariant on the subspace $\mbox{ Im }p$ and $g_0=p-(l')^2$. By the equality on $g_i^{-1}$ we can see that on $\mbox{ Im }p$ $g_0^{-1}-g_1^{-1}=id$ and on $\mbox{ Ker }p$ $g_0^{-1}=g_1^{-1}=-id$. $g_0$ is a regular transformation of $\mbox{ Im }p$ thus the proof of Lemma 9 shows that there exist an invertible selfadjoint transformation $\widetilde{l}$ on $\mbox{ Im }p$ to $\mbox{ Im }p$ for which $g_0=id|_{\mbox{ Im }p}-\widetilde{l}^2$. Extract this transformation to an  $l:V\longrightarrow V$ transformation by the equalities:
$$
l(u)=\left\{\begin{array}{cc}
\widetilde{l}(u) & \mbox{ if } u\in \mbox{ Im }p\\
u & \mbox{ if } u\in \mbox{ Ker }p
\end{array}\right.
$$
Now for an element $u$ of $\mbox{ Im }p$
$$
g_0(u)=u-\widetilde{l}^2(u)=u-{l}^2(u)=(p-(l')^2)(u),
$$
and for a $v\in \mbox{ Ker }p$
$$
g_0(v)=-v=(p-(lp+(id-p))^2)(v)=(p-(l')^2(v),
$$
showing that
$g_0=p-(l')^{2}$ and $\Phi_0$ defined by the equality
$$
\langle p(x),p(x)\rangle  -\langle l'(x),l'(x)\rangle  =0.
$$
It remains to show that indeed $l'(\Phi_0)=\Phi_1$. But by Lemma 7' $l'(\Phi_0)$ defined by a transformation $\overline{g}_1$ with the property
$
g_0^{-1}-(\overline{g}_1)^{-1}=p
$ and thus  $(\overline{g}_1)^{-1}=g_1^{-1}$.
\qed

Now we have

\begin{statement}[Lemma 10']
By the notation of Lemma 9' there is $\delta >  0$ such that $p-\lambda g_0=(l')_{\lambda} ^{2}$ which smoothly depends on $\lambda $, for $-\delta< \lambda< 1+\delta$.
\end{statement}

The easy extraction of the proof of Lemma 10 we omit. Lemma 11 says about the definition of a smooth family of regular transformation corresponding to two "confocal quadrics which are of the same type". Our method of generalization leads to the following analogous statement:

\begin{statement}[Lemma 11']
Suppose that $p,\Phi_0,\Phi_1,g_0,g_1,l$ are as in Lemma 9' and its proof. Then there is a smooth family $l_\lambda$ of transformations with $l_0=id$ and $l_1=l$, such that the quadric $\Phi_\lambda=l'(\Phi_0)$ is defined by the transformation $g_\lambda $ with
$$
g_\lambda ^{-1}=g_0^{-1}-\lambda p.
$$
All quadrics $\Phi_\lambda$ are $p$-quadric belonging to the projection pencil of $p$. They restriction to $\mbox{ Im }p$ are orthogonally intersect the path $l_\lambda(u)$ of a point $u\in \Phi_0\cap \mbox{ Im }p$.
\end{statement}

\proof By our definitions $g_0=p-(l')^2$ and $g_1^{-1}=g_0^{-1}-p$. Consider $\lambda g_0$ instead of $g_0$ then $\lambda g_0 =-id $ on $\mbox{ Ker }p$.  We can define $l_{\lambda}$ by
$$
\lambda g_0=id-(l_{\lambda})^2 \mbox{ on } \mbox{ Im }p,
$$
and by
$$
l_\lambda = id|_{\mbox{ Ker }p} \mbox{ on } \mbox{ Ker }p.
$$
By Lemma 10', $(l'_\lambda)=lp+(id-p) $ exists and depends smoothly on $\lambda $. Now $(l'_\lambda )p=p(l'_\lambda )$ because on $\mbox{ Im }p$ it is an identity and on $\mbox{ Ker }p$ the values of both sides are zero, respectively. (We note that $l'_\lambda=id|_{\mbox{ Ker }p}$ is invariant on $\mbox{ Ker }p$ by its definition.) So we have for non-zero $\lambda$ on $\mbox{ Im }p$
$$
(l_\lambda) g_0=(l_\lambda )\lambda^{-1}(id|_{\mbox{ Im }p}-(l_{\lambda})^2)=\lambda^{-1}(id|_{\mbox{ Im }p}-(l_{\lambda})^2)(l_\lambda )=g_0(l_{\lambda}).
$$
So on $\mbox{ Im }p$
$$
g_\lambda ^{-1}=l_\lambda g_0^{-1}l_\lambda =(l_\lambda )^2g_0^{-1}=(id|_{\mbox{ Im }p}-\lambda g_0)g_0^{-1}=g_0^{-1}-\lambda id.
$$
On the other hand for an element of $\mbox{ Ker }p$ by definition
$$
g_\lambda ^{-1}=g_0^{-1}
$$
showing that on $V$ we have
$$
g_\lambda ^{-1}=l'_\lambda g_0^{-1}l'_\lambda =g_0^{-1}-\lambda p
$$
as we stated. From Lemma 7' we can see that $\Phi_0$ and $l'_{\lambda}(\Phi_0)$ gives a projection pencil corresponding to the projection $p$. Finally, we have to prove the statement on orthogonality. We can get the derivative of the mapping $l_\lambda(x):\mathbb{R}\longrightarrow V$  if we use the direct product structure of $V$. Let $w=u+v\in V$ where $u\in \mbox{ Im }p$ and $v\in \mbox{ Ker }p$. Then we have:
$$
\lambda g_0(u+v)=(u-l_{\lambda}(u)l_{\lambda}(u))-v.
$$
Differentiating this relation with respect to $\lambda $ we get
$$
g_0(u+v)=-2\dot{l_{\lambda}}(u)l_{\lambda}(u),
$$
implying
$$
\dot{l_{\lambda}}(u)=-\frac{1}{2}g_0(u+v){(l_{\lambda}(u))}^{-1}=-\frac{1}{2}g_{\lambda }(u+v){l_{\lambda}(u)},
$$
as in the paper \cite{stachel}. So if $v=0$ then the tangent hyperplane of $\Phi_\lambda$ in $l_{\lambda}(u)$ has the gradient vector $g_{\lambda }l_{\lambda}(u)$ implying that in $P(V)$ the corresponding point is conjugate to the tangent hyperplane with respect to the identity quadric of $Im(P)$ namely to $p$.
\qed

A natural variation of Lemma 12 is valid, too. We have to rewrite now the corresponding assumptions and proof.

\begin{statement}[Lemma 12']
By the notation of the previous lemmas, if $\Psi \neq \Phi_0$ are $p$-quadrics belonging to the same projection pencil (of $p$ and ${g_0}^{-1}$) and $u\in \Phi_0\cap \Psi \cap \mbox{ Im }p$, then $l_\lambda (u)\in \Psi$.
\end{statement}

\proof
We have
$$
x\in \Phi_0 \Longleftrightarrow \langle x,g_0(x)\rangle  =0 \mbox{ and } x\in \Psi \Longleftrightarrow \langle x,g_{\mu}(x)\rangle  =0.
$$
By definition $g_{\mu}^{-1}=g_0^{-1}-\mu p$ with $\mu \neq 0$. Consider the following expression:
$$
\lambda g_0g_\mu^{-1}- \mu l_\lambda g_\mu l_\lambda g_\mu^{-1} -(\lambda -\mu )g_\mu g_\mu^{-1}=
$$
$$
=\lambda g_0(g_0^{-1}-\mu p)-\mu l_\lambda g_\mu l_\lambda g_\mu^{-1}-(\lambda -\mu )id=
$$
$$
=\lambda id-\lambda \mu g_0 p-\mu l_\lambda ^2-(\lambda -\mu )id=
$$
$$
=\lambda id-\lambda \mu g_0 p-\mu (p-\lambda g_0)-(\lambda -\mu )id=
$$
$$
\mu (-p+id)(\lambda g_0+id)
$$
For a points of $\mbox{ Im }p$ it is zero since $p=id$. Thus we also have that
$$
\lambda g_0- \mu l_\lambda g_\mu l_\lambda  -(\lambda -\mu )g_\mu =0
$$
on $\mbox{ Im }p$. Thus on $\mbox{ Im }p$ we get
$$
\mu \langle l_\lambda (u),g_\mu l_\lambda (u)\rangle  =\mu \langle u,l_\lambda g_\mu l_\lambda (u)\rangle  =\lambda \langle u,g_0 (u)\rangle  -(\lambda -\mu)\langle u,g_\mu (u)\rangle  =0,
$$
showing that $l_\lambda (u)\subset \Psi\cap \mbox{ Im }p$.
\qed

Now the proof of Theorem 2 of \cite{stachel} changing the applied lemmas to their variations mentioned and proved in this paper we get the proof of our Theorem. In fact, by Lemma 9' (Statement 5), there exists $l$ such that $\Phi_1=l'(\Phi_0)$ with the transformation $l'=lp+(id-p)$ and the equation of $\Phi_0$ is given by
$$
\langle p(x),p(x)\rangle  -\langle l'(x),l'(x)\rangle  =0.
$$ 
Lemma 11' (Statement 7) shows the existence of $\Phi_\lambda $ and $l_\lambda$. By Lemma 6' (Statement 2), $l_\lambda$ has the Ivory property. At least, Lemma 12' (Statement 8) shows the statement about the quadric $\Psi $, if it exists.

\end{document}